\documentclass[12pt]{amsart}
\usepackage{amsmath,amssymb,amsthm}
\makeatletter
\newcommand{\symbitem}[1]{\item[#1]%
\renewcommand{\@currentlabel}{#1}\ignorespaces}
\makeatother

\makeatletter\@addtoreset{equation}{section} \makeatother

\newtheorem{theorem}[equation]{Theorem}
\newtheorem*{theorem*}{Theorem}

\newtheorem{lemma}[equation]{Lemma}
\newtheorem{corollary}[equation]{Corollary}

\theoremstyle{definition}
\newtheorem{definition}[equation]{Definition}

\theoremstyle{remark}
\newtheorem{remark}[equation]{Remark}

\newcommand{\tit}{On rationality of nonsingular threefolds\\ 
with a pencil of Del Pezzo surfaces of degree $4$}
\title{\tit}
\author{C.\,Shramov}
\thanks{The work was partially supported by RFFI grants
$04-01-00613$ and $05-01-00353$.}

\email{shramov@mccme.ru}

\begin{document}

\maketitle

\begin{abstract}
We prove a criterion of nonsingularity of a complete intersection of two fiberwise
quadrics in 
${\mathbb{P}}_{{\mathbb{P}}^1}({\mathcal{O}}(d_1)\oplus\ldots\oplus{\mathcal{O}}(d_5))$. 
As a corollary we derive the following addition to the Alexeev theorem on 
rationality of standard Del Pezzo fibrations of degree $4$ over 
${\mathbb{P}}^1$: 
we prove that any fibration of this kind with the topological Euler 
characteristic $\chi(X)=-4$ is rational.
\end{abstract}

\section{Introduction}

Rationality questions of Del Pezzo fibrations of degree $d$ over 
${\mathbb{P}}^1$ were studied by many authors.
All the varieties of this type with $d\geqslant 5$ are rational. 
Fibrations with $d=1$, $2$ and $3$ were studied, in particular, 
in~\cite{Pukhlikov}, \cite{Grin1}, \cite{Grin2}, \cite{Grin3},
\cite{ChPrSh}, \cite{Cheltsov}, 
and as a result a nearly complete solution of rationality problem of nonsingular
Del Pezzo fibrations was obtained for $d=1$ and (under generality assumptions) 
for $d=2$,~$3$.

In case of degree $4$ there is a well-known statement.

\begin{theorem*}[Alexeev, see~\cite{Alexeev}]
Let $V$ be a standard Del Pezzo fibration\footnote{See section~\ref{setup}
for a definition.} of degree $4$ over
${\mathbb{P}}^1$. If the topological Euler characteristic 
$\chi(V)\neq 0, -8, -4$, then $V$ is nonrational, if $\chi(V)=0, -8$, then 
$V$ is rational; finally, if $\chi(V)=-4$, then $V$ if and only if
its intermediate Jacobian is a Jacobian of a curve.
\end{theorem*}

The main goal of this paper is to prove the following addition to Alexeev
theorem.

\begin{theorem}\label{theorem:Alexeev-theorem-refined}
Any standard Del Pezzo fibration of degree $4$ over ${\mathbb{P}}^1$ 
with topological Euler characteristic $-4$ is rational.
\end{theorem}

This statement will appear as a corollary of the description of intersections
of fiberwise quadrics in the scrolls
with topological Euler characteristic equal to $-4$ 
(see.~Lemma~\ref{lemma:interesting-varieties} after the necessary notations are
introduced), that in turn is a corollary of the nonsingularity criterion for
intersections of fiberwise quadrics in a fivefold scroll 
(Theorem~\ref{maintheorem}).

\smallskip
The author is grateful to I.\,A.\,Cheltsov for his nondecreasing interest 
to this work and for numerous useful advice, to V.\,A.\,Iskovskikh for 
attention to the work and to S.\,S.\,Galkin for valuable discussions.

\section{Definitions, notations and the statements of the main results}
\label{setup}

All varieties are assumed to be defined over the field of complex numbers 
${\mathbb{C}}$.

\begin{definition}[see, e.\,g.,~\cite{Alexeev}]
A nonsingular threefold $X$ is called a standard Del Pezzo fibration over 
$\mathbb{P}^1$, if $X$ is endowed with a structure of Del Pezzo 
fibration over $\mathbb{P}^1$ with normal fibers, and $\rho(X)=2$.
\end{definition}

Let $\phi:X\to\mathbb{P}^1$ be a standard Del Pezzo fibration of degree $d$ 
over $\mathbb{P}^1$. If $d\geqslant 3$, then  
$X$ is naturally embedded into the projectivisation  
${\mathbb{P}}_{\mathbb{P}^1}({\mathcal{E}})$ of the bundle 
${\mathcal{E}}=\phi_*\omega_X^{-1}$.
For example, if $d=4$, then $X$ is embedded into some scroll 
$Y={\mathbb{P}}_{{\mathbb{P}}^1}({\mathcal{O}}(d_1)\oplus\ldots\oplus{\mathcal{O}}(d_5))$
as a nonsingular complete intersection of two divisors, 
such that the restriction of each of them on a fiber is a quadric in 
${\mathbb{P}}^4$. 
Hence we may always assume that  $X$ is embedded into $Y$ 
as a complete intersection of two fiberwise quadrics.
On the other hand, if $X$ is a complete intersection of this kind, and 
$X$ is nonsingular with $\rho(X)=2$, then $X$ is a standard Del Pezzo 
fibration of degree $4$. It means that to describe Del Pezzo fibrations
of degree $4$ over ${\mathbb{P}}^1$ explicitly the first thing we need
is a way to decide if an intersection of two general members
of given linear systems $|D_1|$ and $|D_2|$ of fiberwise quadrics on $Y$ 
is nonsingular (or, equivalently, if there exist such divisors in these 
linear systems that their intersection is nonsingular).

Let us fix some notations to formulate the answer to the latter question.
Let $Y={\mathbb{F}}(d_1, \ldots, d_5)={\mathbb{P}}_{{\mathbb{P}}^1}({\mathcal{O}}(d_1)\oplus\ldots\oplus{\mathcal{O}}(d_5))$,
where $d_1\geqslant\ldots\geqslant d_5=0$ (the definition and the main 
properties of $Y$ are given, for example, in~\cite{Reid}). Let $M_Y$ be a
tautological invertible sheaf on $Y$, $L_Y$ be a fiber of the natural projection
$\varphi:Y\to{\mathbb{P}}^1$
(we'll often denote them by $M$ and $L$ instead of $M_Y$ and $L_Y$ if no 
ambiguity is likely to arise). 
Let $D_i\in |2M+b_iL|$, $i=1, 2$, be general divisors, given by equations
$f_1(x, t)=\sum\alpha_{ij}x_ix_j$ and 
$f_2(x, t)=\sum\beta_{ij}x_ix_j$ respectively, where $x_i$ are standard
coordinates in a fiber of $\varphi$, and 
$\alpha_{ij}=\alpha_{ij}(t_0, t_1)$ 
(resp. $\beta_{ij}(t_0, t_1)$) are the polynomials in the coordinates 
on the base of degree $d_i+d_j+b_1$ (resp. $d_i+d_j+b_2$). 
We are interested in the conditions on $d_i, b_j$ implying that
$X=D_1\cap D_2$ is a nonsingular variety. 

Let $Y_2, \ldots, Y_5$ be negative subscrolls of $Y$
(i.\,e. $Y_i$ is a variety given by the equations  
$x_1=\ldots=x_{i-1}=0$). It's easy to see that the base locus
$Bs|D_1|$ must coincide with $Y_3$, $Y_4$, $Y_5$ or $\varnothing$.
We'll always assume that $b_1\leqslant b_2$ 
(in particular, $Bs|D_2|\subset Bs|D_1|$).

The following theorem gives criterion for $X$ being nonsingular in terms of
the parameters $d_i$, $b_j$. 
Note that our argument is analogous to one applied to the case of Del
Pezzo fibrations of degree $3$ in~\cite[Lemma~26]{Cheltsov}  
and~\cite[Proposition~31, 32]{BrCoZu}.

\begin{theorem}\label{maintheorem}
Under the assumptions listed above a general variety $X$ 
is nonsingular if and only if one (and hence only one) of the following sets
of conditions holds.

\begin{enumerate}
\item\label{emp} $b_1\geqslant 0$.
\item\label{y5} $b_1<0$, $2d_4+b_1\geqslant 0$, and one of the following sets 
of conditions holds
\begin{enumerate}
\item
\label{y5case1} $d_1+b_1<0$, $b_2=0$.
\item
\label{y5case2} $d_1+b_1\geqslant 0$, $b_2\geqslant 0$.
\item
\label{y5case3} $b_1=-d_1$, $b_2<0$, $d_3+b_2\geqslant 0$.
\item
\label{y5case4} $b_1=-d_1$, $b_2=-d_2$, $d_2>d_3$.
\item
\label{y5case5} $d_1+b_1>0$, $d_2+b_1\geqslant 0$, 
$b_2<0$, $d_3+b_2\geqslant 0$.
\end{enumerate}

\item\label{y4} $2d_4+b_1<0$, $2d_3+b_1\geqslant 0$, and one of the following 
sets of conditions holds
\begin{enumerate}
\item
\label{y4case1} $d_1>d_2$, $d_4>0$, $b_1=-(d_1+d_4)$, $b_2=0$.
\item
\label{y4case2} $d_1>d_2$, $d_4=0$, $b_1=-d_1$, $b_2=0$. 
\item
\label{y4case3} $d_1>d_2+d_4$, $d_4>0$, 
$b_1=-d_1$, $b_2=-2d_4$, $d_3+b_2\geqslant 0$ or $d_2+b_2=0$.
\item
\label{y4case4} $d_1+b_1<0$, $d_2+d_4+b_1\geqslant 0$, $b_2=0$.
\item
\label{y4case5} $d_1+b_1\geqslant 0$, $d_2+d_4+b_1\geqslant 0$, $d_2+b_1<0$, 
$d_3+d_4+b_1<0$, $b_2\geqslant 0$.
\item
\label{y4case6} $d_1>d_2$, $d_1>d_3+d_4$, $d_2+d_4>d_1$,
$b_1=-d_1$, $b_2<0$, $2d_4+b_2\geqslant 0$, $d_3+b_2\geqslant 0$ or 
$ d_2+b_2=0$.
\item
\label{y4case7} $d_1=d_2+d_4$, $d_2>d_3$, $d_4>0$, $b_1=-d_1$, $b_2<0$,
$2d_4+b_2\geqslant 0$, $d_3+b_2\geqslant 0$ or $d_2+b_2=0$.
\item
\label{y4case8} $d_1=d_2+d_4$, $d_2=d_3+d_4$, $d_4>0$, 
$b_1=-d_1$, $b_2=-d_2$.
\item
\label{y4case9} $d_1+b_1\geqslant 0$, $d_3+d_4+b_1\geqslant 0$, $d_2+b_1<0$, 
$b_2\geqslant 0$. 
\item
\label{y4case10} $d_3+d_4\geqslant d_1$, $b_1=-d_1$, 
$b_2<0$, $2d_4+b_2\geqslant 0$,
$d_3+b_2\geqslant 0$ или $d_2+b_2=0$.
\item
\label{y4case11} $d_2+b_1\geqslant 0$, $d_3+d_4+b_1\geqslant 0$, $b_2\geqslant 0$. 
\item
\label{y4case12} $d_2+b_1\geqslant 0$, $d_3+d_4+b_1\geqslant 0$, $b_2<0$, 
$2d_4+b_2\geqslant 0$, $d_3+b_2\geqslant 0$.
\item
\label{y4case13} $d_1=d_2=d_3>0$, $d_4=0$, $b_1=b_2=-d_1$.
\end{enumerate}

\item\label{y3} $2d_3+b_1<0$, $2d_2+b_1\geqslant 0$, and one of the 
following sets of conditions holds
\begin{enumerate}
\item
\label{y3case1} $d_1+d_4=d_2+d_3$, $d_3>d_4>0$, $b_1=-(d_1+d_4)$, $b_2=0$.
\item
\label{y3case2} $d_1=d_2>d_3=d_4>0$, $b_1=-(d_1+d_4)$, $b_2=0$.
\item
\label{y3case3} $d_1=d_2+d_3$, $d_3>d_4>0$, $b_1=-d_1$, $b_2=-2d_4$,
$d_3+b_2\geqslant 0$ or $d_2+b_2=0$.
\item
\label{y3case4} $d_1=d_2+d_3$, $d_3>d_4=0$, $b_1=-d_1$, $b_2=0$.
\item
\label{y3case5} $d_1=d_2$, $d_3=d_4=0$, $b_1=-d_1$, $b_2=0$.
\item
\label{y3case6} $d_4>0$, $d_3=d_4$, $d_2=2d_4$, $d_1=3d_4$, $b_1=-3d_4$,
$b_2=-2d_4$.
\end{enumerate}
\end{enumerate}
\end{theorem}

\begin{remark}
It's easy to check that for each set of conditions in Theorem~\ref{maintheorem}
there is a variety $X$ such that those conditions hold for $X$.
\end{remark}

Theorem~\ref{maintheorem} is proved in section~\ref{proofofmaintheorem}. 
Its awkwardness is partially compensated by the following corollary 
(see section~\ref{application}).

\begin{lemma}\label{lemma:interesting-varieties}
Let $X$ be a standard Del Pezzo fibration of degree $4$ over 
${\mathbb{P}}^1$, and $\chi(X)=-4$. Then $X$ is isomorphic to a complete 
intersection of two fiberwise quadrics in a scroll, and there are only the
following possibilities for parameters $d_i$, $b_j$.
\begin{itemize}
\item[$(X_1)$] $d_1=d_2=d_3=d_4=0$, $b_1=0$, $b_2=1$ 
(case~\ref{emp} of Theorem~\ref{maintheorem}).
\item[$(X_2)$] $d_1=2$, $d_2=d_3=d_4=1$, $b_1=-2$, $b_2=-1$
(case~\ref{y5case3} of Theorem~\ref{maintheorem}).
\end{itemize}
\end{lemma}

Finally, Theorem~\ref{theorem:Alexeev-theorem-refined} is an implication
of Lemma~\ref{lemma:interesting-varieties} provided that the varieties 
$X_1$ and  $X_2$ are rational (this is checked in section~\ref{application}).

\section{Preliminaries}

\begin{lemma}
The variety $X$ is nonsingular if and only if the following conditions hold:
\begin{itemize}
\symbitem{$(*)$}\label{singcond} 
the intersection of $Bs|D_1|\setminus Bs|D_2|$ with $D_2\cap Sing D_1$ is empty,
\symbitem{$(**)$}\label{commonbscond}  
in any point of $Bs|D_2|$ 
the vectors ${\mathrm{grad}_x}(f_1)$ and ${\mathrm{grad}_x}(f_2)$ are not 
proportional.
\end{itemize}
\end{lemma}
\begin{proof}
If $X$ is nonsingular, than the conditions \ref{singcond} and 
\ref{commonbscond} apparently hold.
Assume that the conditions \ref{singcond} and \ref{commonbscond} hold. 
In the points of
$D_1\cap (D_2\setminus Bs|D_2|)$ the variety $D_1$ is nonsingular by
\ref{singcond}, and the divisor $D_2$ is movable, so $X$ is nonsingular 
outside $Bs|D_2|$ by Bertini theorem. In the points of $Bs|D_2|$
the variety $X$ is nonsingular if and only if the vectors
${\mathrm{grad}}(f_1)=({\mathrm{grad}_x}(f_1), {\mathrm{grad}}_t(f_1))$ and 
${\mathrm{grad}}(f_2)=({\mathrm{grad}_x}(f_2), {\mathrm{grad}}_t(f_2))$ 
are not proportional in any point of $Bs|D_2|$.
Let $D'_i$, $i=1, 2$, be divisors given by equations 
$\frac{\partial f_i}{\partial t_0}=0$. Then $|D'_i|=|2M+(b_i-1)L|$, and
$Bs|D'_i|\supset Bs|D_i|$.
Hence ${\mathrm{grad}}_t(f_1)$ and ${\mathrm{grad}}_t(f_2)$ are zero
on $Bs|D_2|$, and the nonproportionality condition may be rewritten as
\ref{commonbscond}.
\end{proof}

To check that $X$ is nonsingular we shall use the conditions
\ref{singcond} and \ref{commonbscond} everywhere below. 

Let's fix some notations. Let $M_4$ denote the $2\times 3$-matrix
$$M_4=\left(
\begin{array}{ccc}
\alpha_{14}x_4+\alpha_{15}x_5&\alpha_{24}x_4+\alpha_{25}x_5&
\alpha_{34}x_4+\alpha_{35}x_5\\
\beta_{14}x_4+\beta_{15}x_5&\beta_{24}x_4+\beta_{25}x_5&
\beta_{34}x_4+\beta_{35}x_5
\end{array}
\right),$$
and $M_5$ --- the $2\times 4$-matrix
$$M_5=\left(
\begin{array}{cccc}
\alpha_{15}&\alpha_{25}&\alpha_{35}&\alpha_{45}\\
\beta_{15}&\beta_{25}&\beta_{35}&\beta_{45}
\end{array}
\right).$$
Let $m^{(l)}_{ij}$, $1\leqslant i< j\leqslant 3$, $l=4, 5$ denote a 
$2\times 2$-minor of the matrix $M_l$, containing its $i$'th and $j$'th columns.

Consider the following conditions:
\begin{itemize}
\symbitem{$(**)_4$}\label{minorcond4} the polynomials $m^{(4)}_{ij}$ 
have no common zeros on $Y_4$.
\symbitem{$(**)_5$}\label{minorcond5} 
the polynomials $m^{(5)}_{ij}$ have no common zeros on $Y_5$.
\end{itemize}
(They are useful because the condition 
\ref{commonbscond}\ is equivalent to \ref{minorcond4} if $Bs|D_2|=Y_4$, 
and to \ref{minorcond5} if $Bs|D_2|=Y_5$.)

\begin{lemma}\label{lemma:M4-minors}
Under the assumptions made above the condition \ref{minorcond4} 
holds if and only if 
$b_1=-d_1=-(d_2+d_4)$, $b_2=-d_2=-(d_3+d_4)$.
\end{lemma}
\begin{proof}
Clearly, to satisfy \ref{minorcond4} we need the last column of $M_4$ to be 
nonzero, hence it is necessary $d_3+d_4+b_2\geqslant 0$. Moreover, 
if $\beta_{25}=0$ or $\alpha_{15}=0$, then $m^{(4)}_{ij}$ are zero on $Y_5$, 
hence \ref{minorcond4} implies $d_2+b_2\geqslant 0$ and $d_1+b_1\geqslant 0$. 
Finally, if
$\alpha_{24}=0$, then \ref{minorcond4} also doesn't hold, so \ref{minorcond4}\
implies $d_2+d_4+b_1\geqslant 0$.

Now assume that the inequalities listed above hold.
Consider the surface $Y_4\cong{\mathbb{F}}(d_4, 0)$, the divisors 
$C_{ij}\in |2M_{Y_4}+(d_i+d_j+b_1+\nolinebreak b_2)L_{Y_4}|$ on $Y_4$, 
given by the equations 
$m^{(4)}_{ij}=0$, the divisor
$A\in |M_{Y_4}+\nolinebreak (d_1+\nolinebreak b_1)L_{Y_4}|$, 
given by the equation $\alpha_{14}x_4+\alpha_{15}x_5=0$, 
and the divisor $B\in |M_{Y_4}+\nolinebreak (d_1+\nolinebreak b_2)L_{Y_4}|$,
given by the equation $\beta_{14}x_4+\beta_{15}x_5=0$. 
The divisors $C_{12}$, $A$ and $B$ are movable, $C_{12}$ intersects $C_{13}$ 
and $A$ intersects $B$ transversally. Hence, $C_{12}\cap C_{13}$ 
consists of $C_{12}C_{13}=4(d_1+d_4+b_1+b_2)+2(d_2+d_3)$ points, and $A\cap B$ 
consists of $AB=2d_1+d_4+b_1+b_2$ points. 

If $C_{12}C_{13}=0$, then the condition
\ref{minorcond4} holds; the equality  
$$0=C_{12}C_{13}=4(d_1+b_1)+2(d_2+b_2)+2(d_3+d_4+b_2)+2d_4$$ 
is equivalent (under the above assumptions) to 
$b_1=b_2=-d_1=-d_2=-d_3$, $d_4=0$. 
Let $C_{12}C_{13}\neq 0$. Note that the vanishing of $m^{(4)}_{12}$ and 
$m^{(4)}_{13}$ implies the vanishing of $m^{(4)}_{23}$ 
provided that the first column of $M_4$ is nonzero.
Hence if $C_{12}C_{13}\neq 0$, then for the condition \ref{minorcond4} to hold
it is necessary and sufficient that $C_{12}\cap C_{13}\subset A\cap B$. 
Since $A\cap B\subset C_{12}\cap C_{13}$, the latter is equivalent to
\begin{multline*} 
0=(4(d_1+d_4+b_1+b_2)+2(d_2+d_3))-(2d_1+d_4+b_1+b_2)=\\
=2(d_1+d_2+d_3)+3(d_4+b_1+b_2)=\\
=2(d_1+b_1)+(d_2+b_2)+(d_2+d_4+b_1)+2(d_3+d_4+b_2).
\end{multline*}
Under the above assumptions this equality is equivalent to 
$b_1=-d_1=-(d_2+d_4)$, 
$b_2=-d_2=-(d_3+d_4)$. To get the final answer note that these conditions 
are weaker than those that appeared in the case $C_{12}C_{13}=0$.
\end{proof}

\begin{lemma}\label{lemma:M5-minors}
Under the assumptions made above the condition \ref{minorcond5} holds 
if and only if either the conditions $d_1+b_1\geqslant 0$ and 
$d_3+b_2\geqslant 0$ hold, or the conditions
$d_1+b_1=0$, $d_2+b_2=0$ hold together with at least one of the conditions 
$d_1+b_2=0$ and $d_2+b_1<0$.
\end{lemma}
\begin{proof}
If $d_1+b_1<0$, then the first row of $M_5$ is zero, so $d_1+b_1\geqslant 0$
is necessary for \ref{minorcond5}. If $d_1+b_1\geqslant 0$, then 
\ref{minorcond5} holds either when there are two nonzero minors $m^{(5)}_{ij}$ 
(we may assume that these are $m^{(5)}_{12}$ and $m^{(5)}_{13}$), 
i.\,e. when $\deg(\beta_{25})=d_3+b_2\geqslant 0$; 
or when
$\deg(m^{(5)}_{12})=\deg(\alpha_{15}\beta_{25}-\alpha_{25}\beta_{15})=0$,
i.\,e. when $\deg(\alpha_{15})=d_1+b_1=0$, $\deg(\beta_{25})=d_2+b_2=0$ 
and one of the following conditions holds: $\deg(\beta_{15})=d_1+b_2=0$ or 
$\deg(\alpha_{25})=d_2+b_1<0$.
\end{proof}

The following statement will simplify the calculations 
in section~\ref{proof:Y3}.

\begin{lemma}\label{lemma:Y3-sing}
Let $Bs|D_1|=Y_3$, $d_2+d_3+b_1\geqslant 0$. Then for the condition 
\ref{singcond} to hold it is necessary that $2(d_1+d_2+d_3+d_4)+4b_1+b_2=0$. 
If the latter holds and either $d_1+d_4+b_1\geqslant 0$, $b_2\geqslant 0$, or 
$d_1+b_1\geqslant 0$,
$2d_4+b_2\geqslant 0$, then this condition is also sufficient.
\end{lemma}
\begin{proof}
If $Bs|D_1|=Y_3$, then all the components of ${\mathrm{grad}_x}(f_1)$ with the possible exception of
$\frac{\partial f_1}{\partial x_1}=\alpha_{13}x_3+\alpha_{14}x_4+\alpha_{15}x_5$ and
$\frac{\partial f_1}{\partial x_2}=\alpha_{23}x_3+\alpha_{24}x_4+\alpha_{25}x_5$, are zero on $Y_3$. Since $d_2+d_3+b_1\geqslant 0$, the polynomials 
$\alpha_{23}$ and $\alpha_{13}$ are nonzero.

Let $A_1$ and $A_2$ be the divisors given by the equations 
$\alpha_{13}x_3+\alpha_{14}x_4+\alpha_{15}x_5=0$
and $\alpha_{23}x_3+\alpha_{24}x_4+\alpha_{25}x_5=0$ respectively. Since
${\mathrm{Sing}}(D_1)=A_1\cap A_2\cap Y_3$, the condition \ref{singcond} 
means that  
$D_2\cap A_1\cap A_2\cap Y_3=\varnothing$. 
In particular,  $A_1A_2D_2\vert_{Y_3}=0$ must hold, i.\,e., as 
$|A_1|=|M_{Y_3}+(d_1+b_1)L_{Y_3}|$, $|A_2|=|M_{Y_3}+(d_2+b_1)L_{Y_3}|$, and 
$Y_3\cong{\mathbb{F}}(d_3, d_4, 0)$, the latter condition may be rewritten as
$2(d_1+d_2+d_3+d_4)+4b_1+b_2=0$. 

If $d_1+d_4+b_1\geqslant 0$, then $A_1$ and $A_2$ have no common components, 
and their intersection is an effective curve on $Y_3$, 
hence if $Y_5\not\subset Bs|D_2|$ the condition $A_1A_2D_2\vert_{Y_3}=0$ 
implies $D_2\cap A_1\cap A_2\cap Y_3=\varnothing$. 
If $d_1+b_1\geqslant 0$, then $A_1$ and $A_2$ have no common components, 
and $Y_5\not\subset A_1$, so that if $Y_4\not\subset Bs|D_2|$ the condition
$A_1A_2D_2\vert_{Y_3}=0$ implies $D_2\cap A_1\cap A_2\cap Y_3=\varnothing$.
\end{proof}

\section{Proof of Theorem~\ref{maintheorem}}\label{proofofmaintheorem}

\subsection{Case $Bs|D_1|=\varnothing$}

This case occurs if and only if $b_1\geqslant 0$. The conditions 
\ref{singcond}  and \ref{commonbscond}  hold automatically
(case~\ref{emp} of Theorem~\ref{maintheorem}).

\subsection{Case $Bs|D_1|=Y_5$}

This case occurs if and only if $b_1<0$, $2d_4+b_1\geqslant 0$.
All the components of ${\mathrm{grad}_x}(f_1)$ with the possible exception 
of 
$\frac{\partial f_1}{\partial x_1}=\alpha_{15}$, $\frac{\partial f_1}{\partial x_2}=\alpha_{25}$,
$\frac{\partial f_1}{\partial x_3}=\alpha_{35}$ and 
$\frac{\partial f_1}{\partial x_4}=\alpha_{45}$ are zero on $Y_5$. 

Let us consider several possibilities.

\subsubsection{Case $d_1+b_1<0$.} 
We have ${\mathrm{grad}_x}(f_1)=0$ on $Y_5$. The condition \ref{commonbscond} 
implies
$Bs|D_2|=\varnothing$, and \ref{singcond}  implies that 
$D_2\cap Y_5=\varnothing$, i.\,e. $b_2=0$
(case~\ref{y5case1} of Theorem~\ref{maintheorem}). 
This condition is apparently sufficient for \ref{singcond}  and 
\ref{commonbscond}. 

\subsubsection{Case $d_1+b_1=0$.} 
Under this assumption we have $\deg(\alpha_{15})=0$. If 
$Bs|D_2|=\varnothing$ (i.\,e. $b_2\geqslant 0$), then \ref{singcond}  and
\ref{commonbscond} hold automatically
(case~\ref{y5case2} of Theorem~\ref{maintheorem}). 
If $Bs|D_2|=Y_5$ ($b_2<0$), it is necessary and sufficient to check the 
condition \ref{commonbscond}  on $Y_5$. Since any pair of polynomials  
$\alpha_{ij}$, $\beta_{kl}$ has no common zeros on $Y_5$, the condition 
\ref{commonbscond} is equivalent to \ref{minorcond5}.
Hence by Lemma~\ref{lemma:M5-minors} the condition \ref{commonbscond} 
holds either if $d_3+b_2\geqslant 0$ 
(case~\ref{y5case3} of Theorem~\ref{maintheorem}), 
or if $d_2+b_2=0$ and one of the following conditions holds: 
$d_1+b_2=0$ or $d_2+b_1<0$;
it's easy to see that under the assumptions made above 
$d_1+b_2=0$ is equivalent to $b_1=b_2$, and $d_2+b_1<0$ is equivalent to 
$b_1<b_2$, hence one of these two holds automatically --- this gives
case~\ref{y5case4} of Theorem~\ref{maintheorem}. 

\subsubsection{Case $d_1+b_1>0$.} 
If $d_2+b_1<0$, then $D_1$ has isolated singularities on $Y_5$, and by 
\ref{singcond} it is necessary that $Bs|D_2|=\varnothing$,
i.\,e. $b_2\geqslant 0$. The condition $b_2\geqslant 0$ is apparently  
sufficient for \ref{singcond} and \ref{commonbscond}\
(regardless to $d_2+b_1<0$). 
This is case~\ref{y5case2} of Theorem~\ref{maintheorem}.
If $d_2+b_1\geqslant 0$ and $b_2<0$, then $Bs|D_2|=Y_5$, 
and \ref{commonbscond} is equivalent to \ref{minorcond5},
i.\,e. (since $d_1+b_1\neq 0$) by Lemma~\ref{lemma:M5-minors} 
the condition \ref{commonbscond} is equivalent to $d_3+b_2\geqslant 0$
(case~\ref{y5case5} of Theorem~\ref{maintheorem}).

\subsection{Case $Bs|D_1|=Y_4$}

This case occurs if and only if $2d_4+b_1< 0$, $2d_3+b_1\geqslant 0$.
Under this assumption all the components of ${\mathrm{grad}_x}(f_1)$
with the possible exception of
$\frac{\partial f_1}{\partial x_1}=\alpha_{14}x_4+\alpha_{15}x_5$, 
$\frac{\partial f_1}{\partial x_2}=\alpha_{24}x_4+\alpha_{25}x_5$ and 
$\frac{\partial f_1}{\partial x_3}=\alpha_{34}x_4+\alpha_{35}x_5$
are zero on $Y_4$.

Let us consider several possibilities.

\subsubsection{Case $d_1+d_4+b_1<0$.} 
We have ${\mathrm{Sing}}(D_1)=Y_4$; in particular, ${\mathrm{Sing}}(D_1)$ 
contains a line $l$ in a fiber of $\varphi$. Since $D_2\in |2M+b_2L|$, 
we have $l\cap D_2\neq\varnothing$, and the condition \ref{singcond} 
doesn't hold.

\subsubsection{Case $d_1+d_4+b_1\geqslant 0$, $d_1+b_1<0$, $d_2+d_4+b_1<0$.} 
We have $\alpha_{15}=\alpha_{24}=\alpha_{25}=\alpha_{34}=\alpha_{35}=0$,
and if $\deg(\alpha_{14})>0$, then the divisor $D_1$ is singular along some 
line in a fiber of $\varphi$, so \ref{singcond} implies  
$\deg(\alpha_{14})=0$, i.\,e. $d_1+d_4+b_1=0$. Furthermore, as $D_1$ is singular
along $Y_5$, the condition \ref{singcond} implies $b_2=0$.
These conditions are apparently sufficient for \ref{singcond} and 
\ref{commonbscond}\
(case~\ref{y4case1} of Theorem~\ref{maintheorem}). 

\subsubsection{Case $d_1+b_1\geqslant 0$, $d_2+d_4+b_1<0$.} 
We have $\alpha_{24}=\alpha_{25}=\alpha_{34}=\alpha_{35}=0$, 
and ${\mathrm{Sing}}(D_1)=C$, where the curve $C$ is given by the equation 
$\alpha_{14}x_4+\alpha_{15}x_5=0$ on $Y_4$. Since $C\neq\varnothing$, 
the condition \ref{singcond} implies $Bs|D_2|\subset Y_5$, i.\,e.
$2d_4+b_2\geqslant 0$. The condition \ref{singcond} also implies 
$$0=D_2C=2(d_1+d_4+b_1)+b_2=2(d_1+b_1)+2d_4+b_2,$$ 
i.\,e. $d_1+b_1=0$, $2d_4+b_2=0$. Since $C\not\supset Y_5$, this is sufficient
for \ref{singcond}. If $d_4=0$, then \ref{commonbscond} holds automatically
(case~\ref{y4case2} of Theorem~\ref{maintheorem}), 
and if $d_4>0$, then $Bs|D_2|=Y_5$, 
and by Lemma~\ref{lemma:M5-minors} the condition \ref{commonbscond} 
holds either if $d_3+b_2\geqslant 0$, or if $d_2+b_2=0$ 
(case~\ref{y4case3} of Theorem~\ref{maintheorem}). 

\subsubsection{Case $d_1+b_1<0$, $d_2+d_4+b_1\geqslant 0$.} 
We have $\alpha_{15}=\alpha_{25}=\alpha_{35}=0$. Since the polynomials 
$\alpha_{14}$ and $\alpha_{24}$ have no common zeros,  
${\mathrm{Sing}}(D_1)=Y_5$. 
The condition \ref{singcond} implies $b_2=0$.
The latter is apparently sufficient for \ref{singcond} and \ref{commonbscond}\
(case~\ref{y4case4} of Theorem~\ref{maintheorem}). 

\subsubsection{Case $d_1+b_1\geqslant 0$, $d_2+d_4+b_1\geqslant 0$, 
$d_2+b_1<0$, $d_3+d_4+b_1<0$.}
We have $\alpha_{25}=\alpha_{34}=\alpha_{35}=0$, and $D_1$ has at most isolated 
singularities. If $d_1+b_1>0$, then $D_1$ has some singular points 
on $Y_5$, hence \ref{singcond} implies $b_2\geqslant 0$.
The latter is apparently sufficient for \ref{singcond} and \ref{commonbscond} 
regardless to the assumption $d_1+b_1>0$
(case~\ref{y4case5} of Theorem~\ref{maintheorem}). 

If $d_1+b_1=0$, $b_2<0$ and $d_2+d_4+b_1>0$, then $D_1$ has some singular 
points on
$Y_4\setminus Y_5$, and by \ref{singcond} it is necessary that 
$Bs|D_2|\neq Y_4$,
i.\,e. $2d_4+b_2\geqslant 0$, hence $Bs|D_2|=Y_5$. 
By Lemma~\ref{lemma:M5-minors} in this case the condition \ref{commonbscond} 
holds if and only if either $d_3+b_2\geqslant 0$, or $d_2+b_2=0$
(case~\ref{y4case6} of Theorem~\ref{maintheorem}). 

If $d_1+b_1=0$, $b_2<0$ and $d_2+d_4+b_1=0$, then $D_1$ is nonsingular. 
If under these assumptions $Bs|D_2|=Y_5$, i.\,e. $2d_4+b_2\geqslant 0$, then
by Lemma~\ref{lemma:M5-minors} 
the condition \ref{commonbscond} holds if and only if 
either $d_3+b_2\geqslant 0$, or $d_2+b_2=0$
(case~\ref{y4case7} of Theorem~\ref{maintheorem}). 
If $Bs|D_2|=Y_4$,
i.\,e. $2d_4+b_2<0$, then \ref{commonbscond} is equivalent to \ref{minorcond4}.
By Lemma~\ref{lemma:M4-minors} 
this is equivalent to $b_1=-d_1=-(d_2+d_4)$, $b_2=-d_2=-(d_3+d_4)$
(case~\ref{y4case8} of Theorem~\ref{maintheorem}).

\subsubsection{Case $d_1+b_1\geqslant 0$, $d_3+d_4+b_1\geqslant 0$, 
$d_2+b_1<0$.} 
We have $\alpha_{25}=\alpha_{35}=0$. If $d_1+b_1>0$, then $D_1$ has 
(isolated) singularities on $Y_5$, hence \ref{singcond} implies
$b_2\geqslant 0$. The latter is sufficient for  
\ref{singcond} and \ref{commonbscond} regardless to the assumption $d_1+b_1>0$
(case~\ref{y4case9} of Theorem~\ref{maintheorem}).

If $d_1+b_1=0$, $b_2<0$, then $D_1$ is nonsingular. If $Bs|D_2|=Y_5$, i.\,e.
$2d_4+b_2\geqslant 0$, then by Lemma~\ref{lemma:M5-minors} 
the condition \ref{commonbscond} holds if and only if either 
$d_3+b_2\geqslant 0$, or $d_2+b_2=0$
(case~\ref{y4case10} of Theorem~\ref{maintheorem}). 
If $Bs|D_2|=Y_4$, then \ref{commonbscond} is equivalent to \ref{minorcond4}.
By Lemma~\ref{lemma:M4-minors} the latter implies $d_2=d_3+d_4$,
contradicting the assumptions $d_2+b_1<0$, $d_3+d_4+b_1\geqslant 0$. 

\subsubsection{Case $d_3+d_4+b_1\geqslant 0$, $d_2+b_1\geqslant 0$.} 
The divisor $D_1$ is nonsingular. If $b_2\geqslant 0$, then \ref{commonbscond} 
holds automatically
(case~\ref{y4case11} of Theorem~\ref{maintheorem}). If $b_2<0$, 
$2d_4+b_2\geqslant 0$ (i.\,e. $Bs|D_2|=Y_5$), then, since the latter implies 
$b_2>b_1$, by Lemma~\ref{lemma:M5-minors} the condition \ref{commonbscond} 
holds only if $d_3+b_2\geqslant 0$
(case~\ref{y4case12} of Theorem~\ref{maintheorem}). 
If $2d_4+b_2<0$ (i.\,e. $Bs|D_2|=Y_4$), 
the condition \ref{commonbscond} is equivalent to \ref{minorcond4}. 
By Lemma~\ref{lemma:M4-minors} the latter is equivalent to 
$b_1=-d_1=-(d_2+d_4)$, $b_2=-d_2=-(d_3+d_4)$. 
Under the current assumptions this implies $b_1=b_2$, $d_1=d_2=d_3>0$, 
$d_4=0$ (case~\ref{y4case13} of Theorem~\ref{maintheorem}).

\subsection{Case $Bs|D_1|=Y_3$}\label{proof:Y3}

This case occurs if and only if $2d_3+b_1<0$, $2d_2+b_1\geqslant 0$.
Under this assumption all the components of ${\mathrm{grad}_x}(f_1)$
with the possible exception of 
$\frac{\partial f_1}{\partial x_1}=\alpha_{13}x_3+\alpha_{14}x_4+\alpha_{15}x_5$ and
$\frac{\partial f_1}{\partial x_2}=\alpha_{23}x_3+\alpha_{24}x_4+\alpha_{25}x_5$are zero on $Y_3$.

Let us consider several possibilities.

\subsubsection{Case $d_2+d_3+b_1<0$.} 
We have $\alpha_{23}=\alpha_{24}=\alpha_{25}=0$,
and ${\mathrm{Sing}}(D_1)$ contains a line in a general fiber of $\varphi$, 
hence \ref{singcond} doesn't hold.

\subsubsection{Case $d_2+d_3+b_1\geqslant 0$, $d_2+d_4+b_1<0$, $d_1+b_1<0$.} 
We have $\alpha_{15}=\alpha_{24}=\alpha_{25}=0$. If 
$\deg(\alpha_{23})>0$, then ${\mathrm{Sing}}(D_1)$ contains a curve in a fiber
of $\varphi$ over a zero of the polynomial $\alpha_{23}$. Hence
\ref{singcond} implies $d_2+d_3+b_1=0$. In this case 
${\mathrm{Sing}}(D_1)\subset Y_4$.
If $\deg(\alpha_{14})\neq 0$ (in particular, if this degree is negative, i.\,e.
the polynomial $\alpha_{14}$ is zero), then ${\mathrm{Sing}}(D_1)$ contains
a line in a fiber of $\varphi$, hence it is necessary that $d_1+d_4+b_1=0$. 
In this case  
${\mathrm{Sing}}(D_1)=Y_5$, and for \ref{singcond} and \ref{commonbscond} 
to hold the equality $b_2=0$ is necessary and sufficient 
(case~\ref{y3case1} of Theorem~\ref{maintheorem}).

\subsubsection{Case $d_2+d_4+b_1\geqslant 0$, $d_1+b_1<0$.} 
We have $\alpha_{15}=\alpha_{25}=0$, and $Y_5\subset{\mathrm{Sing}}(D_1)$. 
Hence \ref{singcond} implies $b_2=0$. By Lemma~\ref{lemma:Y3-sing} 
the condition \ref{singcond} is equivalent to
$$0=2(d_1+d_2+d_3+d_4)+4b_1+b_2=2(d_2+d_4+b_1)+2(d_1+d_3+b_1),$$
i.\,e. $d_2+d_4+b_1=d_1+d_3+b_1=0$. Since $b_2=0$, the condition
\ref{commonbscond} holds automatically (case~\ref{y3case2} 
of Theorem~\ref{maintheorem}).

\subsubsection{Case $d_2+d_3+b_1\geqslant 0$, $d_2+d_4+b_1<0$, $d_1+b_1\geqslant 0$.} 
We have $\alpha_{24}=\alpha_{25}=0$. If $d_1+b_1>0$, then 
${\mathrm{Sing}}(D_1)\cap Y_5\neq\varnothing$, and the inequality 
$b_2\geqslant 0$ must hold. 
By Lemma~\ref{lemma:Y3-sing} the condition \ref{singcond} is equivalent to
$$0=2(d_1+d_2+d_3+d_4)+4b_1+b_2=2(d_1+b_1)+2(d_2+d_3+b_1)+2d_4+b_2>0,$$
a contradiction. Hence it is necessary that $d_1+b_1=0$. Since $d_2+d_3+b_1>0$
implies that ${\mathrm{Sing}}(D_1)$ contains a line in a fiber of $\varphi$, 
the equality $d_2+d_3+b_1=0$ is also necessary.

By Lemma~\ref{lemma:Y3-sing} the condition \ref{singcond} implies
$$0=2(d_1+d_2+d_3+d_4)+4b_1+b_2=2(d_1+b_1)+2(d_2+d_3+b_1)+2d_4+b_2,$$
that is equivalent to $2d_4+b_2=0$, i.\,e. (again by Lemma~\ref{lemma:Y3-sing})
under the assumptions made above \ref{singcond} is equivalent to $2d_4+b_2=0$.
If $d_4>0$, then $Bs|D_2|=Y_5$, and \ref{commonbscond} is equivalent to 
\ref{minorcond5}, and as $d_2+b_1<0$ and $d_1+b_1=0$ the latter condition 
holds if and only if $d_3+b_2\geqslant 0$ or
$d_2+b_2=0$ by Lemma~\ref{lemma:M5-minors} 
(case~\ref{y3case3} of Theorem~\ref{maintheorem}). 
If $d_2=0$, i.\,e. $Bs|D_2|=\varnothing$, we have $b_2=0$, which apparently
implies \ref{commonbscond} (case~\ref{y3case4} 
of Theorem~\ref{maintheorem}).

\subsubsection{Case $d_2+d_4+b_1\geqslant 0$, $d_1+b_1\geqslant 0$.} 
Since ${\mathrm{Sing}}(D_1)\neq\varnothing$,
we have $Bs|D_2|\neq Y_3$, i.\,e. $2d_3+b_2\geqslant 0$. 
By Lemma~\ref{lemma:Y3-sing} the condition \ref{singcond} implies 
$$0=2(d_1+d_2+d_3+d_4)+4b_1+b_2=2(d_1+b_1)+2(d_2+d_4+b_1)+2d_3+b_2,$$
i.\,e. $d_2+d_4+b_1=d_1+b_1=2d_3+b_2=0$. If under these assumptions $d_3>d_4$, 
then $Bs|D_2|=Y_4$, and by Lemma~\ref{lemma:M4-minors} the condition 
\ref{commonbscond} implies $b_2=-(d_3+d_4)$, i.\,e. $d_3=d_4$, a contradiction.
Hence we have $d_3=d_4$, and by Lemma~\ref{lemma:Y3-sing} the assumptions made 
above are sufficient for \ref{singcond}. 

The condition $2d_4+b_2=0$ means that $Bs|D_2|\subset Y_5$. 
If $Bs|D_2|=\varnothing$, then $b_2\geqslant 0$ implies $b_2=d_3=d_4=0$;
in this case the condition \ref{commonbscond} holds automatically
(case~\ref{y3case5} of Theorem~\ref{maintheorem}). If $Bs|D_2|=Y_5$,
i.\,e. $d_4>0$, then \ref{commonbscond} is equivalent to \ref{minorcond5}. 
Since $d_2+b_1<0$ and
$d_3+b_2<0$, the condition \ref{minorcond5} is equivalent to $d_2+b_2=0$ 
by Lemma~\ref{lemma:M5-minors} (case~\ref{y3case6} 
of Theorem~\ref{maintheorem}).

\section{Приложение к расслоением на поверхности Дель Пеццо степени $4$}
\label{application}

Applying Theorem~\ref{maintheorem}, we get the following
\begin{corollary}\label{corollary:chi-4}
Let the assumptions made in section~\ref{setup} hold and let 
$X$ be nonsingular with topological Euler characteristic $\chi(X)=-4$.
The one of the following hold.
\begin{itemize}
\item[($X_1$)] $d_1=d_2=d_3=d_4=0$, $b_1=0$, $b_2=1$ 
(case~\ref{emp} of Theorem~\ref{maintheorem}).
\item[($X_2$)] $d_1=2$, $d_2=d_3=d_4=1$, $b_1=-2$, $b_2=-1$
(case~\ref{y5case3} of Theorem~\ref{maintheorem}).
\item[($X_3$)] $d_1=4$, $d_2=3$, $d_3=2$, $d_4=1$, $b_1=-4$, $b_2=-3$
(case~\ref{y4case8} of Theorem~\ref{maintheorem}).
\end{itemize}
\end{corollary}
\begin{proof}
It's easy to check that $\chi(X)=-16\sum d_i-20b_1-20b_2+16$ 
(for example, it is implied by~\cite[Example~3.2.11]{Fulton}). 
All that remains is to solve the equation 
$-16\sum d_i-20b_1-20b_2+16=-4$ together with the sets of equations and inequalities on $d_i$, $b_j$ from Theorem~\ref{maintheorem}.
\end{proof}

Hence all the varieties we are interested in are contained in three families
listed above, and it is sufficient to check rationality of a general member 
of each of these families.

\begin{remark}\label{remark-on-X3-Picard}
In case ($X_3$) both $Bs|D_1|$ and $Bs|D_2|$ contain $Y_4$.
Consider the generic fiber of $X_3$ as a surface over the field 
${\mathbb{C}}(t_0)$ of rational functions on a line. This surface contains 
a line defined over ${\mathbb{C}}(t_0)$, hence the relative Picard number 
$\rho(X_3/{\mathbb{P}}^1)$ is not less than $2$.
\end{remark}

Lemma~\ref{lemma:interesting-varieties} is immediately implied by
Corollary~\ref{corollary:chi-4} and 
Remark~\ref{remark-on-X3-Picard}.

\begin{lemma}\label{lemma:rationality}
Varieties $X_1$ and $X_2$ are rational.
\end{lemma}
\begin{proof}
The projection of $X_1\subset{\mathbb{P}}^4\times{\mathbb{P}}^1$ in 
${\mathbb{P}}^4$ gives a morphism onto a threefold quadric 
$Q\subset{\mathbb{P}}^4$ that is birational since in this case $D_2\in |2M+L|$;
rationality of $X_1$ follows immediately. 

To prove that $X_2\subset Y={\mathbb{F}}(2, 1, 1, 1, 0)$ is rational consider
the projection $\pi:Y\to Y'={\mathbb{F}}(2, 1, 1, 1)$ from the curve 
$Y_5\subset Y$. It gives a birational map of $X_2$ onto a divisor 
$X'_2\in |3M_{Y'}-3L_{Y'}|$. Identifying $Y'$ with 
$Y''={\mathbb{F}}(1, 0, 0, 0)$ we identify the variety 
$X'_2\subset Y''$ with a divisor from the linear system $|3M_{Y''}|$ on $Y''$.
Finally, the contraction $\sigma:Y''\to{\mathbb{P}}^4$ of the negative subscroll
$Y''_2$ represents $Y''$ as a blowup of a plane  
${\mathbb{P}}^2\subset{\mathbb{P}}^4$, and gives a birational morphism of 
$X'_2$ onto a cubic in ${\mathbb{P}}^4$. Hence to prove the rationality 
of $X_2$ it suffices to check that this cubic is singular. 

Let us show that the singularities arise already on $X'_2$. To see that we'll
find all the curves contracted by the map $\pi$. They are lines
in the fibers of $\varphi$, passing through $Y_5$. Each line $l$ with this 
property, written down in the coordinates  
$(x_1, \ldots, x_4)$ in the affine space $x_5=1$, is given parametrically
by $(x_1, \ldots, x_4)=s(v_1, v_2, v_3, v_4)$. 
The condition $l\subset D_1\cap D_2$ means that for any $s$ we have
\begin{multline*}
((\alpha_{11}x_1+\alpha_{12}x_2+\alpha_{13}x_3+\alpha_{14}x_4)x_1+\\
+(\alpha_{22}x_2^2+\alpha_{33}x_3^2+\alpha_{44}x_4^2+
\alpha_{23}x_2x_3+\alpha_{24}x_2x_4+\alpha_{34}x_3x_4))s^2+\\
\shoveright{+\alpha_{15}x_1s=0,}\\
\shoveleft{((\beta_{11}x_1+\beta_{12}x_2+\beta_{13}x_3+\beta_{14}x_4)x_1+}\\
+(\beta_{22}x_2^2+\beta_{33}x_3^2+\beta_{44}x_4^2+
\beta_{23}x_2x_3+\beta_{24}x_2x_4+\beta_{34}x_3x_4))s^2+\\
+(\beta_{15}x_1+(\beta_{25}x_2+\beta_{35}x_3+\beta_{45}x_4))s=0.
\end{multline*}

Since $\alpha_{15}$ is a nonzero constants, these equalities may be rewritten as
\begin{gather}
x_1=0,\label{x1-equation}\\
\alpha_{22}x_2^2+\alpha_{33}x_3^2+\alpha_{44}x_4^2+
\alpha_{23}x_2x_3+\alpha_{24}x_2x_4+\alpha_{34}x_3x_4=0,\label{q-alpha}\\
\beta_{22}x_2^2+\beta_{33}x_3^2+\beta_{44}x_4^2+
\beta_{23}x_2x_3+\beta_{24}x_2x_4+\beta_{34}x_3x_4=0,\label{q-beta}\\
\beta_{25}x_2+\beta_{35}x_3+\beta_{45}x_4=0.\label{l-beta}
\end{gather}
In particular, nothing is contracted in a general fiber. 
On the other hand, there are two fibers containing one contracted line each:
the equations~\ref{x1-equation}, \ref{q-alpha} 
and~\ref{l-beta} (that do not depend on $t$) give two (up to proportionality)
possible values of a vector $v=(v_1, v_2, v_3, v_4)$,
and for each of them the equation~\ref{q-beta} (that is linear in $t$) 
gives exactly one value of $t$, such that a line corresponding to the vector $v$
in the fiber over $t$ is contained in $D_1\cap D_2$. Hence the variety
$X'_2$ has two singular points (that are simple double points). 
\end{proof}

\begin{remark}
Rationality of the variety $X_3$ is also easy to prove --- $X_3$ contains
a surface $Y_4$, and the projection from $Y_4$ maps $X_3$ birationally on a 
rational variety ${\mathbb{F}}(4, 3, 2)$.
\end{remark}

Lemma~\ref{lemma:interesting-varieties} and Lemma~\ref{lemma:rationality}
prove Theorem~\ref{theorem:Alexeev-theorem-refined}.

\end{document}